\documentclass[12pt,twoside]{article}
\usepackage{amsmath,amssymb,mathrsfs,graphicx}
\usepackage[small,sc]{caption} 

\newtheorem{theorem}{Theorem}[section]
\newtheorem{lemma}[theorem]{Lemma}
\newtheorem{proposition}[theorem]{Proposition}
\newtheorem{corollary}[theorem]{Corollary}

\newtheorem{rmrk}[theorem]{Remark}

\DeclareMathAlphabet{\mathbfit}{OML}{cmm}{b}{it}

\makeatletter
\@addtoreset{equation}{section}
\makeatother

\newenvironment{remark}
{\begin{rmrk} \em}
{\end{rmrk}}

\newcommand{\me} {measure}

\newcommand{\erg} {ergodic}
\newcommand{\sy} {system}

\newcommand{\dsy} {dynamical system}

\newcommand{\R} {\mathbb{R}}

\newcommand{\Z} {\mathbb{Z}}
\newcommand{\N} {\mathbb{N}}
\newcommand{\qed} {\hfill {\small Q.E.D.} \par\medskip}
\newcommand{\skippar} {\par\medskip}
\newcommand{\ds} {\displaystyle}
\newcommand{\proof} {\noindent \textsc{Proof.} }
\newcommand{\proofof}[1] {\noindent \textsc{Proof of {#1}.} }
\newcommand{\article}[3] {\textsc{{#1}}, {\itshape {#2}}, {{#3}}.}
\newcommand{\book}[3] {\textsc{{#1}}, {\itshape {#2}}, {{#3}}.}
\newcommand{\vol} {\textbf}

\newcommand{\rset}[2] {\left\{ #1 \: \left| \: #2 \right. \! \right\} }

\newcommand{\symmdiff} {\triangle}

\newcommand{\into} {\longrightarrow}

\renewcommand{\iff} {if and only if\ }

\newcommand{\newfig}[4] {
\bigskip
\begin{figure}[htbp]
  \centering
  \includegraphics[width=#3]{#2}
  \begin{minipage}[t]{0.80\linewidth} 
    \caption{#4}
    \protect\label{#1}
  \end{minipage}
\end{figure}

}

\newcommand{\ps} {\mathcal{M}}   % phase space
\newcommand{\sca} {\mathscr{A}}   % generic sigma-algebra
\newcommand{\scb} {\mathscr{B}} 
\newcommand{\scn} {\mathscr{N}}   % null sigma-algebra
\newcommand{\scp} {\mathscr{P}}   % generic partition
\newcommand{\scq} {\mathscr{Q}}   % partition for the invertible case
\newcommand{\sct} {\mathscr{T}}   % tail sigma-algebra
\newcommand{\fs} {\mathcal{I}}   % fiber set
\newcommand{\fa} {\Psi}   % fiber action
\newcommand{\ec} {\mathcal{E}}   % ergodic component

\begin{document}

\title{\textbf{Extensions of exact and K-mixing dynamical systems}}

\author{
\scshape
Daniele Galli\thanks{
Dipartimento di Matematica, Universit\`a di Bologna,
Piazza di Porta San Donato 5, 40126 Bologna, Italy. 
E-mails: \texttt{daniele.galli7@unibo.it}, \texttt{marco.lenci@unibo.it}.}
\ and
Marco Lenci$^*$\thanks{
Istituto Nazionale di Fisica Nucleare,
Sezione di Bologna, Via Irnerio 46,
40126 Bologna, Italy.}
}

\date{February 2025%
\footnote{Post-publication. Fixed minor imprecisions in version at
  J.\ Stat.\ Phys.\ \vol{190} (2023), no.~1, Paper No. 21.}
}

%%% MODIFICATIONS OVER THE PUBLISHED VERSION
% Clarified proof of Lem. 2.1
% Thms 2.3 + 2.5, (iv) + (vi): added assumption that \mu is invariant
% Prop. 4.1: added assumption that \psi is \scb_o-measurable
% modified Rmk 4.2 as previously incorrect

\maketitle

\begin{abstract}
  We consider extensions of non-singular maps which are exact, respectively 
  K-mixing, or at least have a decomposition into positive-measure exact, 
  respectively K-mixing, components. The fibers of the extension spaces have 
  countable (finite or infinite) cardinality and the action on them is assumed 
  surjective or bijective. We call these systems, respectively, 
  \emph{fiber-surjective} and \emph{fiber-bijective} extensions. Technically, they 
  are skew products, though the point of view we take here is not the one
  generally associated with skew products. Our main results are an Exact and a 
  K-mixing Decomposition Theorem. The latter can be used to show that a large 
  number of periodic Lorentz gases (the term denoting here general group 
  extensions of Sinai billiards, including Lorentz tubes and slabs, in any 
  dimension) are K-mixing.
  
  \bigskip\noindent 
  \textbf{Mathematics Subject Classification (2020):} 37A40, 37A35;
  37C83, 37D25, 37A20.
  
  \bigskip\noindent
  \textbf{Keywords:} Extensions of dynamical systems, exactness, K-property, 
  decomposition theorem, skew products, Lorentz gas, Lorentz tube, Sinai
  billiard.
  
%  \bigskip\noindent
%  \textbf{Data availability statement:} Data sharing not applicable to this article 
%  as no datasets were generated or analyzed during the current study.
\end{abstract}

\section{Introduction}
\label{sec-intro}

A non-invertible map $T: \ps \into \ps$, which is two-sided non-singular for the 
$\sigma$-finite measure space $(\ps, \sca, \mu)$, is called \emph{exact} if 
\begin{equation}
  \bigcap_{n \in \N} T^{-n} \sca = \scn,
\end{equation}
where $\scn$ is the \emph{null}, or \emph{trivial} $\sigma$-algebra on 
$\ps$, containing only the zero-measure sets and their complements. Its
significance is that the \dsy\ $(\ps, \sca, \mu, T)$ carries vanishing 
information about its initial conditions into the far future; equivalently, it 
retains vanishing information about the state of the system in the remote 
past. In fact, if an observer makes a finite-precision measurement of the 
state of the \sy\ at present time and $\mathscr{O}$ is the (necessarily finite) 
$\sigma$-algebra representing the resolution of the measurement, the 
resolution that the observer obtains on the state of the \sy\ at time
$-n$ is given by $T^{-n} \mathscr{O} \subseteq T^{-n} \sca$, which becomes 
trivially coarse as $n \to \infty$.

If $T$ is an invertible map
the above definition cannot work, in general, because $T^{-n} \sca = \sca$ 
for all $n$. One uses instead the definition of \emph{K-mixing}, a.k.a.\ 
\emph{K-property} (after Kolmogorov), which means 
that there exists a $\sigma$-algebra $\scb \subset \sca$ such that
\begin{equation} \label{def-k}
  T \scb \supset \scb, \qquad \sigma\! \left( \bigcup_{n\in \N} T^n \scb \right) 
  = \sca \mod \mu, \qquad \bigcap_{n \in \N} T^{-n} \scb = \scn,
\end{equation}
where $\sigma( \cdot )$ means the $\sigma$-algebra generated by a family 
of subsets. In this case, the $\sigma$-algebra $\mathscr{O}$ of the 
present-time observation, being finite, can be approximated to any degree 
by a $\sigma$-algebra $\mathscr{O}' \subseteq T^N \scb$, for some $N$. 
This implies that $\bigcap_n T^{-n} \mathscr{O}' = \scn$, leading to the same
interpretation as before.

In this paper we study the exact, respectively K-mixing, components of
extensions of exact, respectively K-mixing, \dsy s, or at least \dsy s whose 
exact/K-mixing components have positive measure. If $T_o$ is an 
endomorphism of the Lebesgue space $(\ps_0, \sca_o, \mu_o)$, we 
consider the \sy\ $(\ps, \sca, \mu, T)$, where $\ps := \ps_o \times \fs$ and 
$\fs$ is a countable (finite or infinite) set. The map $T: \ps \into \ps$ acts as
\begin{equation}
  T(x,i) := (T_o(x), \fa(x,i)),
\end{equation}
where the functions $\fa(x, \cdot): \fs \into \fs$ are surjective or bijective, 
depending on the sought result (of course the two assumptions coincide 
when $\fs$ is finite). As for $\sca$ and $\mu$, they are the natural lifts of 
$\sca_o$ and $\mu_o$, respectively, to $\ps_o \times \fs$ (see Section 
\ref{sec-setup} for details). If the functions $\fa(x, \cdot)$ are surjective, 
respectively bijective, we say that $T$ is a \emph{fiber-surjective}, 
respectively \emph{fiber-bijective}, extension of $T_o$. These two classes 
of \dsy s are quite general. For example, every \emph{group extension}, 
that is, a case where $\fs$ is an abelian group and $\fa(x,i) = i + \psi(x)$, for 
some $\psi: \ps_o \into \fs$, is fiber-bijective. Notable examples are periodic 
Lorentz gases, which are $\Z^d$-extensions of Sinai billiards. The recent 
paper \cite{l4} studies the exact components of certain Markov maps of $\R$,
including $\Z$-extensions of expanding circle maps.

A term that is often used for systems similar to fiber-surjective  maps is 
\emph{skew products}, but when one speaks of a skew product, one usually 
means the case where $\fs$ is an uncountable space with some structure
(measure-theoretic, topological, differentiable, etc.)\ and $T_o$ is an 
automorphism of a measure space (typically a probability space). More 
importantly, the focus tends to be on the action of the maps 
$\fa(T_o^{n-1}(x), \cdot) \circ \cdots \circ \fa(T_o(x), \cdot) \circ 
\fa(x, \cdot)$ on $\fs$, rather than the dynamics of $T$ on the 
whole $\ps_o \times \fs$, as is the case here. (See \cite{k, bk} for 
results on the exactness of certain families of skew products.)

We describe our main results in some detail. If 
$(\ps_0, \sca_o, \mu_o, T_o)$ is exact and $T$ is fiber-surjective, there are 
at most countably many \erg\ components of $T$. Each of them is 
partitioned into $m \in \Z^+ \cup \{ \infty \}$ positive-measure sets, which we 
call \emph{atoms}, and $T$ maps atoms onto atoms. If $m \in \Z^+$, $T$ 
acts as an $m$-cycle between the atoms and $T^m$ is exact on each of them. 
If $m = \infty$, $T$ acts as a shift between the atoms. A completely analogous 
theorem is proved for the case where $(\ps_0, \sca_o, \mu_o, T_o)$ is 
K-mixing and $T$ is fiber-bijective (subject to a natural technical condition).
In this case, for $m \in \Z^+$, $T^m$ is K-mixing on the atoms of an 
$m$-cycle. Theorems of this kind are sometimes called Spectral 
Decomposition Theorems \cite{p, ks}, a phrase we choose to avoid in the 
context of infinite ergodic theory. Here we refer to them, respectively, as the 
Exact and the K-mixing Decomposition Theorem. In addition, assuming that
$T$ is fiber-bijective (even for the case of $T_o$ exact), we are able to provide
extra information on the structure of the \erg\ components and of the
atoms within cycles. This implies in particular that, if $T$ is conservative (i.e.,
recurrent), it is isomorphic to another map $T_1$ on $(\ps, \sca, \mu)$ whose 
atoms are made up of sets $\ps_o \times \{i\}$.

Finally, as an application of the K-mixing Decomposition Theorem, we show 
that, rather generally, a periodic Lorentz gas is K-mixing \iff it is conservative. 
In particular, this establishes the K-property for a large number of periodic 
Lorentz \emph{tubes} and \emph{slabs} (see Section \ref{sec-lg}).

\skippar

The paper is organized as follows. In Section \ref{sec-setup} we present
the mathematical setup, establish the notation that is used throughout 
the paper, and state our main results, which are proved in Section 
\ref{sec-proofs}. In Section \ref{sec-lg}, which is designed to be read 
independently of Section \ref{sec-proofs}, we show how to apply the 
K-mixing Decomposition Theorem to periodic Lorentz gases.

\paragraph{Acknowledgments.} This research was partially supported by the 
PRIN Grant 2017S35EHN, MUR, Italy. It is also part of the authors' activity 
within the UMI Group \emph{DinAmicI} and M.L.'s activity within the Gruppo 
Nazionale di Fisica Matema\-tica, INdAM.

\section{Setup and results}
\label{sec-setup}

Given a \dsy\ $(\ps_o, \sca_o, \mu_o, T_o)$, with $(\ps_o, \sca_o, \mu_o)$ 
a Lebesgue \me\ space and $T_o: \ps_o \into \ps_o$ bimeasurable and 
two-sided non-singular (i.e., $\mu_o(A)=0$ $\Leftrightarrow$ 
$\mu_o(T_o^{-1} A)=0$), we consider its extension $(\ps, \sca, \mu, T)$, 
where:
\begin{itemize}
\item $\ps := \ps_o \times \fs$, for some countable (finite or infinite) $\fs$;
\item $\sca$ is the lift of $\sca_o$ to $\ps$, i.e., $\sca := \sigma ( 
  \rset{A \times \{i\} }{A \in \sca_o,\, i\in \fs} )$;
\item $\mu$ is the lift of $\mu_o$ to $\ps$, i.e., $\mu$ is the \me\ on $\sca$ 
  uniquely defined by $\mu( A \times \{i\} ) = \mu_o(A)$, for all $A \in \sca_o$ 
  and $i\in \fs$;
\item $T$ is a bimeasurable self-map of $\ps$ given by $T(x,i) = 
  (T_o(x), \fa(x,i))$. Depending on the context, we assume that, for 
  $\mu_o$-a.a.\ $x \in \ps_o$,
  \begin{equation} \label{main-hyp}
    \fa_x := \fa(x, \cdot) : \fs \into \fs \ \mbox{is surjective or bijective}.
  \end{equation}
  In the first case we call $T$ a {\bfseries fiber-surjective extension} of $T_o$; 
  in the second case, we call it a {\bfseries fiber-bijective extension} of $T_o$. 
  Clearly, the two conditions are equivalent if $\# \fs < \infty$.
\end{itemize}

Observe that, since $T_o$ is bimeasurable, the bimeasurability of $T$ is 
equivalent to the condition that the level sets of $\fa( \cdot, i)$ are 
measurable, that is, for all $i,j \in \fs$, $\rset{x\in\ps_o} {\fa(x,i)=j} \in \sca_o$.
This is in turn equivalent to the measurability and bimeasurability of 
$\fa$ w.r.t.\ the power set of $\fs$.

We denote by $\pi: \ps \into \ps_o$ the natural projection $\pi(x,i) := x$. By 
construction, $\pi \circ T = T_o \circ \pi$. The main purpose of assumption
(\ref{main-hyp}) is to also have $\pi \circ T^{-1} = T_o^{-1} \circ \pi$, in the
sense of a relation between sets, as follows.

\begin{lemma} \label{lem-comm2}
  If $T$ is a fiber-surjective extension of $T_o$ then, for all $A \in \sca$, 
  \begin{displaymath}
    \pi \, T^{-1} A = T_o^{-1} \, \pi A  \mod \mu.
  \end{displaymath}
\end{lemma}

\proof It is enough to show that $\pi(T^{-1}(y,j)) = T_o^{-1}(\pi(y,j)) = 
T_o^{-1}(y)$, for $\mu$-a.e.\ $(y,j) \in \ps$. Apart from a $\mu_o$-null
set of $x \in \ps_o$, $x \in T_o^{-1}(y)$ implies that there exists $i \in \fs$ 
such that $T(x,i) = (y,j)$, by the fiber-surjectivity of $T$. By non-singularity,
that relation holds for $\mu_o$-a.e.\ $y \in \ps_o$ and for all $j \in \fs$. Also, 
clearly, no $x \not\in T_o^{-1}(y)$ can be such that $T(x,i) = (y,j)$ for some 
$i$. Thus, $\mu$-a.s.\ in $(y,j)$, 
\begin{equation}
  \pi(T^{-1}(y,j)) = \pi( \rset{ (x,i) \in \ps} {T(x,i) = (y,j)} ) = T_o^{-1}(y),
\end{equation}
as claimed.
\qed

\begin{proposition}
  If $T$ is fiber-surjective, $T$ is two-sided non-singular \iff $T_o$ is. 
  If $T$ is fiber-bijective, $T$ preserves $\mu$ \iff $T_o$ preserves $\mu_o$.
\end{proposition}

\proof The two-sided non-singularity is an immediate consequence of 
Lemma \ref{lem-comm2}, since, for all $A \in \sca$, $\mu(A) = 0$ 
$\Leftrightarrow$ $\mu_o(\pi A) = 0$.

As for the second statement, it suffices to check the preservation of $\mu$
on all sets of the type $B \times \{j\}$. By the injectivity of a.a.\ $\fa_x$, cf.\
(\ref{main-hyp}), the sets
\begin{equation}
  B_i := \rset{x \in \ps_o} {T(x,i) \in B \times \{j\} } 
\end{equation}
are disjoint mod $\mu_o$ in $\ps_o$. In other words, $\pi$ is injective mod
$\mu$ on $T^{-1}(B \times \{j\})$, giving 
\begin{equation}
  \mu( T^{-1}(B \times \{j\}) ) = \sum_{i \in \fs} \mu(B_i \times \{i\}) = \mu_o 
  \! \left( \bigsqcup_{i \in \fs} B_i \right) = \mu_o ( T_o^{-1} B ),
\end{equation}
where the last equality follows from the surjectivity of a.a.\ $\fa_x$. Since
$\mu(B \times \{j\}) = \mu_o(B)$, this concludes the proof of the second 
statement.
\qed

\bigskip\noindent
\textbf{Convention.} From now on, all statements about sets are
intended modulo null sets, w.r.t.\ the relevant measure. For example, if
$A,B \in \ps$, the equality $A=B$ means $\mu(A \symmdiff B) = 0$; the 
inclusion $A \subseteq B$ means $\mu(A\setminus B) = 0$; with 
$A \subset B$ implying in addition that $\mu(B\setminus A) > 0$. 
The convention includes $\sigma$-algebras. For instance, $\sca 
\subseteq \scb$ is intended in the sense that the inclusion holds for
the respective completions.

\bigskip

We now present our decomposition theorems and corollaries. In order to do 
so, we recall that, for any non-singular \dsy\ $(\ps, \sca, \mu, T)$, $\ps$ can 
be decomposed into a \emph{conservative part} $\mathcal{C}$ and a 
\emph{dissipative part} $\mathcal{D}$ \cite[\S1.1]{a}. $\mathcal{C}$ is the 
part where Poincar\'e recurrence holds, 
$\mathcal{D}$ is the measurable union of all \emph{wandering sets}. This
is called the \emph{Hopf decomposition} of $T$ and is such that $T^{-1} 
\mathcal{C} \supseteq \mathcal{C}$ and $T^{-1} \mathcal{D} \subseteq 
\mathcal{D}$. It $T$ is measure-preserving, these inclusions become
equalities. It is easy to see that $T^m$, $m\ge 2$, has the same Hopf 
decomposition as $T$. The proofs of all the following results are given in 
Section \ref{sec-proofs}.

\begin{theorem} \label{thm-main-ex}
  Let $(\ps_o, \sca_o, \mu_o, T_o)$ be exact with $(\ps, \sca, \mu, T)$ a 
  fiber-surjective extension, as defined earlier. Then the tail $\sigma$-algebra 
  \begin{displaymath}
    \sct(T) := \bigcap_{n\in\N} T^{-n} \sca
  \end{displaymath}
  is \emph{atomic}, in the sense that it is generated by a partition $\scp$ of 
  $\ps$, with the following properties:
  \begin{itemize}
  \item[(i)] Every $P \in \scp$ (henceforth called an \emph{atom} of $\scp$ 
    or $\sct(T)$) has positive (possibly infinite) \me. In particular, $\scp$ is 
    countable.
    
  \item[(ii)] For all $P \in \scp$ and $n \in \Z$, $T^n P \in \scp$.
   
  \item[(iii)] All \erg\ components are of the form $\ec = \ec_P := 
    \bigcup_{n \in \Z} T^n P$, for some $P \in \scp$. An \erg\ component can 
    comprise a finite number $m$ of atoms, in which case we call it an 
    \emph{$m$-cycle}, or an infinite number of atoms, in which case we call 
    it a \emph{chain}.
     
   \item[(iv)] If $P \in \scp$ and $m \in \Z^+$ are such that $T^m P = P$, 
     then $T^m|_P : P \into P$ is exact. In particular, if $\mu$ is 
     invariant and $\mu(P) < \infty$, then $\ec_P$ belongs to the conservative 
     part of the \sy.
  \end{itemize}
  If $T$ is also fiber-bijective then:
  \begin{itemize}
  \item[(v)] Given an \erg\ component $\ec$, $\# \rset{i \in \fs} 
    {(x,i) \in \ec}$ is constant for $\mu_o$-a.e.\ $x \in \ps_o$. 
    Denoting it $N_\ec \in \Z^+ \cup \{\infty\}$, it follows that 
    $\pi |_\ec: \ec \into \ps_o$ is $N_\ec$-to-1 and onto, whence $\mu(\ec) 
    = N_\ec \, \mu_o(\ps_o)$.
  
  \item[(vi)] If $\ec = \ec_P$ is an $m$-cycle, then also $\# \rset{i \in \fs} 
    {(x,i) \in P}$ is constant for $\mu_o$-a.e.\ $x \in \ps_o$. 
    Denoting it $N_P \in \Z^+ \cup \{\infty\}$, it follows that 
    $\pi |_P: P \into \ps_o$ is $N_P$-to-1 and onto, whence $\mu(P) = 
    N_P \, \mu_o(\ps_o)$. Furthermore, $N_P = N_{P'}$ for all atoms $P'$ 
    of $\ec_P$. If $\mu$ is invariant, this implies $\mu(P) = \mu(P')$ and
    $\mu(\ec_P) = m \mu(P)$.
  
  \item[(vii)] If $T$ is conservative, $T$ is isomorphic to another 
    fiber-bijective extension $T_1$, defined again on $(\ps, \sca, \mu)$, the 
    atoms of whose $\sigma$-algebra are made up of entire \emph{levels} 
    $\ps_o \times \{i\}$.
  \end{itemize} 
\end{theorem}

One need not assume the exactness of the base \sy\ to ensure that a 
fiber-bijective extension has an atomic tail $\sigma$-algebra. If the base 
\sy\ itself has an atomic tail $\sigma$-algebra, this property carries over to 
the extension.

\begin{corollary} \label{cor-main-ex}
  Assertions (i)-(iv) of Theorem \ref{thm-main-ex} also hold in the case 
  where the tail $\sigma$-algebra of $(\ps_o, \sca_o, \mu_o, T_o)$ is atomic, 
  i.e., $\sct(T_o)$ is generated by a partition with positive-\me\ elements. 
  Moreover, if $P$ is an atom of $\sct(T)$, then $\pi P$ is an atom of 
  $\sct(T_o)$.
\end{corollary}

Analogous results hold for the case where $T_o$ is K-mixing. The familiar
reader will see how the extra assumption on the $\scb$-measurability of 
$\fa( \cdot, i)$ is obvious if $T$ is to exploit the K-property of $T_o$.

\begin{theorem} \label{thm-main-k}
  Let $(\ps_o, \sca_o, \mu_o, T_o)$ be K-mixing w.r.t.\ the $\sigma$-algebra
  $\scb_o \subset \sca_o$, with $(\ps, \sca, \mu, T)$ a fiber-bijective 
  extension, as defined earlier. Denote by $\scb$ the lift of $\scb_o$ to $\ps$. 
  Assume that $\fa$ is measurable w.r.t.\ $\scb$ (equivalently,
  $\fa( \cdot, i): \ps_o \into \fs$ is $\scb_o$-measurable for all $i$; 
  equivalently, the level sets of $\fa( \cdot, i )$ are $\scb_o$-measurable).
  Then the $\sigma$-algebra 
  \begin{displaymath}
    \sct_\scb(T) := \bigcap_{n\in\N} T^{-n} \scb
  \end{displaymath}
  is atomic, being generated by a partition $\scp$ of $\ps$ with the following
  properties:
  \begin{itemize}
    \item[(i)] Every $P \in \scp$ (again called an atom of $\scp$ or 
    $\sct_\scb(T)$) has positive (possibly infinite) \me. In particular, $\scp$ is 
    countable.
    
  \item[(ii)] For all $P \in \scp$ and $n \in \Z$, $T^n P \in \scp$.
   
  \item[(iii)] All \erg\ components are of the form $\ec = \ec_P := 
    \bigcup_{n \in \Z} T^n P$, for some $P \in \scp$. 
    
   \item[(iv)] If $P \in \scp$ and $m \in \Z^+$ are such that $T^m P = P$, then 
     $T^m|_P : P \into P$ is K-mixing w.r.t.\ the $\sigma$-algebra
     $\scb \cap P := \rset{A \cap P} {A \in \scb}$. In particular, if $\mu$ is 
     invariant and $\mu(P) < \infty$, then $\ec_P$ belongs to the conservative 
     part of the \sy.

    \item[(v)] Given an \erg\ component $\ec$, $\# \rset{i \in \fs} 
    {(x,i) \in \ec}$ is constant for $\mu_o$-a.e.\ $x \in \ps_o$. 
    Denoting it $N_\ec \in \Z^+ \cup \{\infty\}$, it follows that 
    $\pi |_\ec: \ec \into \ps_o$ is $N_\ec$-to-1 and onto, whence $\mu(\ec) 
    = N_\ec \, \mu_o(\ps_o)$.
  
  \item[(vi)] If $\ec = \ec_P$ is an $m$-cycle, then also $\# \rset{i \in \fs} 
    {(x,i) \in P}$ is constant for $\mu_o$-a.e.\ $x \in \ps_o$. 
    Denoting it $N_P \in \Z^+ \cup \{\infty\}$, it follows that 
    $\pi |_P: P \into \ps_o$ is $N_P$-to-1 and onto, whence $\mu(P) = 
    N_P \, \mu_o(\ps_o)$. Furthermore, $N_P = N_{P'}$ for all atoms $P'$ 
    of $\ec_P$.  If $\mu$ is invariant, this implies $\mu(P) = \mu(P')$ and
    $\mu(\ec_P) = m \mu(P)$.
      
  \item[(vii)] If $T$ is conservative, there is an isomorphism of \dsy s that takes 
    $T$ to another fiber-bijective extension $T_1$ on $(\ps, \sca, \mu)$, and 
    maps $\sct_\scb(T)$ to a $\sigma$-algebra whose atoms are made up of 
    entire levels $\ps_o \times \{i\}$.
  \end{itemize} 
\end{theorem}

\begin{corollary} \label{cor-main-k}
  Assertions (i)-(iv) above also hold in the case where 
  $(\ps_o, \sca_o, \mu_o, T_o)$ has an atomic K-mixing decomposition relative 
  to $\scb_o$, that is, 
  \begin{displaymath}
    T_o \scb_o \supset \scb_o, \qquad \sigma\! \left( \bigcup_{n\in \N} 
    T_o^n \scb_o \right) = \sca_o \mod \mu, \qquad \sct_{\scb_o}(T_o)
    \mbox{ is atomic.}
  \end{displaymath}
  Moreover, if $P$ is an atom of $\sct_\scb(T)$, then $\pi P$ is an atom of 
  $\sct_{\scb_o}(T_o)$.
\end{corollary}

\section{Proofs}
\label{sec-proofs}

\proofof{Theorem \ref{thm-main-ex}} It is a general simple fact that $A \in 
\sct(T)$ \iff $A = T^{-n} \, T^n A$, for all $n \in \N$. For any such $A$, 
Lemma \ref{lem-comm2} gives
\begin{equation}
  \pi A = \pi \, T^{-n} \, T^n A = T_o^{-n} \, T_o^n \, \pi A,
\end{equation}
and so $\pi A \in \sct(T_o) = \scn_o$, by the exactness of $T_o$, where 
$\scn_o$ is the trivial $\sigma$-algebra of $\ps_o$. This implies that if 
$\mu(A) > 0$ (and so $\mu_o(\pi A) > 0$), then 
\begin{equation} \label{pf10}
  \mu(A) \ge \mu_o(\pi A) = \mu_o(\ps_o).
\end{equation}

Since $(\ps_o, \sca_o, \mu_o)$ is a Lebesgue space and $\fs$ is countable, 
$(\ps, \sca, \mu)$ is also a Lebesgue space. Therefore, there exists a unique 
measurable partition $\scp$ that generates $\sct := \sct(T)$ (the uniqueness is 
intended mod $\mu$, according to our convention) \cite{r}. Let $D$
be the union of all the zero-measure elements of $\scp$: $D \in \sct$
because it is the complement of the (necessarily countable) union of all
the positive-measure elements of $\scp$. 

We claim that $\mu(D)=0$. If not, by construction, $D$ is infinitely
divisible in $\sct$, in the sense that any of its positive-measure subsets can 
be split into two positive-measure subsets belonging to $\sct$. In particular, 
there exist $D_1, D_2 \in \sct$, with $\mu(D_1), \mu(D_2) >0$ and 
$D = D_1 \sqcup D_2$. Recursively, for every $n \ge 2$ and
$i_1, \ldots, i_{n-1}, i_n \in \{1,2\}$, there exists a set 
$D_{i_1, \ldots, i_{n-1}, i_n} \in \sct$ such that
\begin{equation} \label{pf20}
  D_{i_1, \ldots, i_{n-1}, i_n} \subset D_{i_1, \ldots, i_{n-1}}, \qquad 
  \mu( D_{i_1, \ldots, i_{n-1}, i_n} ) > 0.
\end{equation}
By (\ref{pf10}), then, $ \mu( D_{i_1, \ldots, i_{n-1}, i_n} ) \ge \mu_o(\ps_o)$.
For every $\bar\imath = (i_1, \ldots, i_n, \ldots) \in \{1,2\}^{\Z^+}$, define
\begin{equation} \label{pf30}
  D_{\bar\imath} := \bigcap_{n=1}^\infty D_{i_1, \ldots, i_n} \in \sct.
\end{equation}
Clearly, $\mu( D_{\bar\imath} ) \ge \mu_o(\ps_o)$, but the sets 
$\{ D_{\bar\imath} \}$ are uncountably many
and disjoint, in contradiction with the fact that $(\ps, \sca, \mu)$ is 
$\sigma$-finite (being Lebesgue). This proves our claim and assertion
\emph{(i)} of Theorem \ref{thm-main-ex}.

Now, recall that $\sct$ is invariant for both $T$ and $T^{-1}$.
If $P, P' \in \scp$ and $n \in \Z$, the case $0 < \mu( T^n P \cap P' ) < 
\mu(P')$ cannot occur, otherwise $T^n P \cap P'$ would be a strict subset 
of the atom $P'$. Also $T^n P \supset P'$ is excluded, otherwise, since $T$
is two-sided non-singular, $T^{-n} P'$ would be a strict subset of the atom $P$.
(Notice that, if $n>0$, we must also use the identity $T^{-n} T^n P = P$.)
It follows that $T^n P$ is either $P'$ or disjoint from $P'$, mod $\mu$, 
yielding \emph{(ii)}.

Assertion \emph{(iii)} is a simple consequence of the fact that any invariant
set $A = T^{-1} A$ belongs to $\sct$ and thus must be a union of atoms. 

As for \emph{(iv)}, it suffices to observe that, for a general $T$, $\sct(T^m) = 
\sct(T)$ for all $m>0$; and $\sct(T|_A) = \sct(T) \cap A$, for any invariant set 
$A$. Under our hypotheses, this implies that
\begin{equation} 
  \sct(T^m|_P) = \sct(T^m) \cap P = \sct(T) \cap P = \scn \cap P,
\end{equation}
the trivial $\sigma$-algebra on $P$. Hence $T^m|_P$ is exact. For the 
second assertion of \emph{(iv)}, we first observe that, by the Poincar\'e 
Recurrence Theorem applied to $T^m |_{\ec_P}$, $P$ belongs to the 
conservative part of $T^m$, which coincides with $\mathcal{C}$, the 
conservative part of $T$. Secondly, any other atom $P'$ of $\ec_P$ can be 
rewritten as $P' = T^n P$, for some $n>0$. Since 
$T \mathcal{C} \subseteq \mathcal{C}$, we have that
$P' \subseteq \mathcal{C}$, ending the proof of \emph{(iv)}.

For the rest of the proof we assume that $T$ is fiber-bijective. For 
$x \in \ps_o$, set $f_\ec (x) := \# \rset{i \in \fs} {(x,i) \in \ec}$. Since $\ec$ is 
$T$-invariant and $T$ is fiber-bijective, $f_{\ec}$ is a $T_o$-invariant 
observable of $\ps_o$, taking values in $\N \cup \{ \infty \}$. By the 
ergodicity of $T_o$, $f_{\ec}$ is a.e.\ equal to a constant value $N_\ec$, 
which cannot be zero because $\mu(\ec)>0$. This proves \emph{(v)}. 

The first part of \emph{(vi)} is proved by the above argument applied to $P$, 
$T^m$ and $T_o^m$. The fact that $N_P = N_{P'}$, whenever $P'$ is an atom
 of $\ec_P$, comes from the fiber-bijectivity of $T^n$, for $n \in \Z^+$ such that
 $T^n P= P'$.

Lastly, we prove \emph{(vii)}. We will construct an isomorphism of measure 
spaces $\phi : \ps \into \ps$ by defining it separately on each atom $P \in \scp$. 
So fix any such $P$ and consider $N_P \in \Z^+ \cup \{\infty\}$ as given by 
\emph{(vi)}. For $x \in \ps_o$, set $P_x := \rset{i \in \fs} {(x,i) \in P}$. Since 
$\fs$ is countable, we endow it with a well-order whereby, for $\mu_o$-a.e.\ 
$x$, $P_x$ can be enumerated as $\{ i_j(x) \}_{j=1}^{N_P}$, where 
$i_j: \ps_o \into \fs$ is recursively defined by 
\begin{align}
  i_1(x) &:= \min \rset{i \in \fs} {(x,i) \in P}, \\
  i_{j+1}(x) &:= \min \rset{i \in \fs} {i > i_j(x), \, (x,i) \in P}.
\end{align}
These expressions prove that the functions $i_j$ are measurable. Now, if 
$N_P \in \Z^+$, set $\Z_{N_P} := \{ 1, 2, \ldots, N_P \}$; if $N_P = \infty$, 
set $\Z_{N_P} := \Z^+$. For $x$ such that $\# P_x = N_P$, set
\begin{equation}
  \phi_P(x,i_j(x)) := (x,j).
\end{equation}
This defines a bijection $\phi_P : P \into \ps_o \times \Z_{N_P}$ mod $\mu$. 
The expression $\phi_P^{-1} (x,j) = (x,i_j(x))$ shows that $\phi_P^{-1} $ is 
measurable, whereas the fact that, for all $B \in \sca_o$ and $j \in \Z_{N_P}$,
\begin{equation}
  \phi_P^{-1}(B \times \{j\}) := \rset{(x,i_j(x))} {x \in B} = \bigsqcup_{k \in \fs}
  (B \cap i_j^{-1}(k)) \times \{k\}.
\end{equation}
proves that $\phi_P$ is measurable and carries the measure $\mu|_P$ to 
the measure $\mu_o \times \#$ on $\ps_o \times\Z_{N_P}$, where $\#$ is 
the counting measure.

We aggregate all bijections $\phi_P$ to construct a bijection $\phi: \ps \into
\bigsqcup_{P \in \scp} \ps_o \times \Z_{N_P}$ (all bijections are intended
mod $\mu$ and $\mu_o \times \#$, respectively). We do so, with a slight 
abuse of notation on the codomain, by defining $\phi(x,i) := \phi_P(x,i)$ 
whenever $(x,i) \in P$. Evidently, $\sum_{P \in \scp} N_P = \# \fs$ and this
cardinality can be finite or countably infinite. In either case, there exists 
a bijection $\gamma: \bigsqcup_{P \in \scp} \Z_{N_P} \into \fs$.
Defining $\Phi: \ps \into \ps$ via $\Phi := (\mathrm{id}, \gamma) \circ \phi$ 
results in an automorphism of the measure space $(\ps, \sca, \mu)$ which,
by construction, acts on the factor $\ps_o$ as the identity. In other words,
it commutes with the projection $\pi$. Setting $T_1 := \Phi \circ T \circ \Phi^{-1}$
gives the sought isomorphism, because $T_1 \circ \pi = \pi \circ T_1$ and 
the atoms of the tail $\sigma$-algebra of $T_1$ are the sets
$\Phi P = \ps_o \times \gamma \!\left( \Z_{N_P} \right)$, for all $P \in \scp$.
\qed

\medskip

\proofof{Corollary \ref{cor-main-ex}} This corollary is in fact a porism, in that
its proof derives not from the statement of Theorem \ref{thm-main-ex}, but from 
its proof, which we refer the reader to.

Let $D$ be the union of all the null elements of $\sct := \sct(T)$. Proving the
absurdity of the claim $\mu(D)>0$ is slightly more involved that in the proof of
Theorem \ref{thm-main-ex}, because $A \in \sct$ no longer implies (\ref{pf10}). 
On the other hand, for any atom $P_o$ of $\sct(T_o)$, $P_o \times \fs \in
\sct$ (because $T^{-n} T^n (P_o \times \fs) =  (T_o^{-n} T_o^n P_o) \times \fs
= P_o \times \fs$). So, for some $P_o$, $D' := D \cap (P_o \times \fs)$ must
have positive measure. But $D'$ is infinitely divisible by construction and, 
for all $A \in \sct$ with $A \subseteq D'$ and $\mu(A)>0$, the inequality 
$\mu(A) \ge \mu_o(\pi A) = \mu_o(P_o)$ holds for the same reasons as for
(\ref{pf10}). This leads to a contradiction as in (\ref{pf20})-(\ref{pf30}). 
Assertion \emph{(i)} is proved. As for \emph{(ii)-(iv)}, the proof of Theorem 
\ref{thm-main-ex} works as is.

Regarding the last assertion of Corollary \ref{cor-main-ex}, we already know
that $\pi P \in \sct(T_o)$. Given an atom $P_o$ of $\sct(T_o)$, it cannot
be $\pi P \cap P_o \subset P_o$, otherwise $P_o$ would not be an atom; and
neither can be that $P_o \subset \pi P$, for in that case $(P_o \times \fs)
\cap P \in \sct(T)$, with $\mu( (P_o \times \fs) \cap P) > 0$. Since $P$ is an 
atom, this would imply $P \subseteq P_o \times \fs$ and thus 
$\pi P \subseteq P_o$, a contradiction. Therefore $\pi P$ must either coincide 
with or be disjoint from $P_o$ (mod $\mu_o$). 
\qed

\medskip

\proofof{Theorem \ref{thm-main-k}} The idea here is to reduce the invertible 
dynamical system $(\ps, \sca, \mu, T)$ to a non-invertible \sy\ whose base 
dynamics is exact, cf.\ \cite[Rmk 3.3]{l3}, and then apply Theorem 
\ref{thm-main-ex}.

Let $\scq_o$ denote the partition of $\ps_o$ associated to $\scb_o$ and 
$[x]$ the element of $\scq_o$ containing a given $x \in \ps_o$. For 
$B \subseteq \ps_o$, denote $[B] := \rset{[x] \in \scq_o} {x \in B}$. Observe 
that this transformation of sets is invertible mod $\mu_o$ on $\scb_o$, in the 
sense that, if $B, B' \in \scb_o$ and $[B] = [B']$, then $B=B'$ mod $\mu_o$ 
(this is one of the properties of the correspondence between 
$\sigma$-algebras and measurable partitions in Lebesgue spaces \cite{r}). 
This shows that the
collection $[\scb_o] := \rset{[B]} {B \in \scb_o}$ is a $\sigma$-algebra for
$\scq_o$ and the measure $\mu_{[\scb_o]}([B]) := \mu_o(B)$ is well-defined 
on it. Finally define $T_{\scb_o}: \scq_o \into \scq_o$ via $T_{\scb_o}([x]) 
= [T_o(x)]$. This definition is well-posed since $T_o$ is K-mixing w.r.t.\
$\scb_o$, implying that $T_o \scq_o$ is a finer partition than $\scq_o$.
Clearly $(\scq_o, [\scb_o], \mu_{[\scb_o]}, T_{\scb_o})$ is a two-sided
non-singular factor of $(\ps_o, \scb_o, \mu_o, T_o)$, which is
non-invertible because $T_o \scq_o$ is strictly finer than $\scq_o$.
The factor property shows that $\sct( T_{\scb_o} ) = [ \sct_{\scb_o} (T_o) ]
= [\scn]$. that is, $T_{\scb_o}$ is exact.

We make the analogous construction for $\scq := \scq_o \times \fs$, which
is evidently the partition of $\ps$ corresponding to $\scb$. We end up with the
system $(\scq, [\scb], \mu_{[\scb]}, T_\scb)$, which is thus a fiber-bijective
$\fs$-extension of $(\scq_o, [\scb_o], \mu_{[\scb_o]}, T_{\scb_o})$. A very
important remark here is that $T_\scb$ is well-defined \iff $\fa$ is 
$\scb$-measurable.

Theorem \ref{thm-main-ex} can be applied to $T_\scb$, implying properties 
\emph{(i)-(vii)} for the dynamics of $T_\scb$ on $\sct(T_\scb)$. But the 
transformation $\scb \ni A \mapsto [A]$ is invertible mod $\mu$, in the sense 
explained earlier, so $P$ is an atom of $\sct_\scb(T)$ \iff $[P]$ is an atom of 
$\sct(T_\scb)$, and they have the same measure, respectively w.r.t.\ $\mu$
and $\mu_{[\scb]}$. This proves all the assertions of Theorem \ref{thm-main-k},
except for the K-mixing of $T^m|_P$ w.r.t.\ $\scb \cap P$, when $P$ is part of
an $m$-cycle.

The triviality of $\sct_{\scb \cap P}(T^m)$ is equivalent to the triviality of 
 $\sct_\scb (T) \cap P$, which is equivalent to the triviality of $\sct(T_\scb)
\cap [P]$, which we have proved. For the two remaining conditions of 
(\ref{def-k}) we state the following lemma, which will be proved at the end
of Section \ref{sec-proofs}.

\begin{lemma} \label{lem-k-mix}
  With the notation $\scb_o \times \fs := \rset{B \times \{i\}} {B \in\scb_o,\, i \in \fs}$,
  the following hold:
  \begin{itemize}
  \item[(i)] For all $n \in \N$, $T^n ( \sigma(\scb_o \times \fs) ) \supseteq 
    \sigma(T_o^n \scb_o \times \fs)$;
  \item[(ii)] $\ds \sigma \!\left( \bigcup_{n \in \N} \sigma(T_o^n\scb_o\times\fs)
    \right) = \sigma \!\left( \sigma \!\left( \bigcup_{n \in \N} T_o^n \scb_o \right)
    \times \fs \right) $.
  \end{itemize}
\end{lemma}

Lemma \ref{lem-k-mix}\emph{(i)}, with $n=1$, gives
\begin{equation} \label{pf50}
  T \scb = T(\sigma(\scb_o \times \fs)) \supseteq \sigma (T_o \scb_o \times \fs) 
  \supset \sigma(\scb_o \times \fs) = \scb,
\end{equation}
which in turn implies $T^m (\scb \cap P) = (T^m \scb) \cap P \supset \scb 
\cap P$. Furthermore,
\begin{equation} \label{pf60}
\begin{split}
  \sigma \!\left( \bigcup_{n\in\N} T^n \scb \right) &= \sigma \! \left(
  \bigcup_{n\in\N} T^n \sigma \!\left( \scb_o \times \fs \right) \right) \\[2pt]
  &\supseteq \sigma \! \left( \bigcup_{n\in\N} \sigma \!\left( T_o^n \scb_o \times 
  \fs \right) \right) \\[2pt]
  &= \sigma \!\left( \sigma \!\left( \bigcup_{n\in\N}T_o^n \scb_o \right)
  \times \fs \right) \\[4pt]
  &= \sigma \!\left( \sca_o \times \fs \right) = \sca,
\end{split}
\end{equation}
where the second line comes from Lemma \ref{lem-k-mix}\emph{(i)}, the third 
line comes from Lemma \ref{lem-k-mix}\emph{(ii)}, and the fourth line comes
from the hypothesis $\sigma \!\left( \bigcup_{n\in\N} T_o^n \scb_o \right) = 
\sca_o$. But (\ref{pf60}) must be a chain of equalities because its leftmost term 
cannot exceed $\sca$. The K-property of $T^m|_P$ is proved.
\qed

\medskip

\proofof{Corollary \ref{cor-main-k}} This proof employs a combination of 
arguments explained in the two preceding proofs. 

One constructs the systems $(\scq_o, [\scb_o], \mu_{[\scb_o]}, T_{\scb_o})$ 
and $(\scq, [\scb], \mu_{[\scb]}, T_\scb)$ as in the proof of Theorem 
\ref{thm-main-k}, and reasons as in the proof of Corollary \ref{cor-main-ex} 
to show that $\sct (T_\scb)$, and therefore $\sct_\scb (T)$, is atomic with
the properties \emph{(i)-(iii)} and part of \emph{(iv)}. The K-property of 
$T^m|_P$, when $P$ is an atom of an $m$-chain, was already established 
in the previous proof. The last claim of Corollary \ref{cor-main-k} is proved 
as in the proof of Corollary \ref{cor-main-ex}.
\qed

\medskip

\proofof{Lemma \ref{lem-k-mix}} We start with inclusion \emph{(i)}. Since 
$T^n$ is invertible, $T^n ( \sigma(\scb_o \times \fs) )$ is a $\sigma$-algebra,
so it suffices to show that $T_o^n \scb_o \times \fs \subseteq T^n( \sigma
( \scb_o \times \fs))$. We may assume $n=1$ (otherwise we use 
$T_o^n, T^n$ in lieu of $T_o, T$). For all $B \in \scb_o$ and $j \in \fs$, we need 
to prove that $T_o B \times \{j\} \in T(\sigma(\scb_o \times \fs))$.
For $i \in \fs$, set $B_i := B \cap \rset{x \in \ps_o} {\fa(x,i) = j}$. Since 
$\fa(\cdot, i)$ is $\scb_o$-measurable, $B_i \in \scb_o$. By fiber-bijectivity, 
$B = \bigsqcup_{i \in \fs} B_i$. Therefore
\begin{equation}
  T_o B \times \{j\} = \left( T_o \bigsqcup_{i \in \fs} B_i \right) \! \times \{j\} =
  T \!\left( \bigsqcup_{i \in \fs} B_i \times \{i\} \right) \in T \!\left( \sigma ( 
  \scb_o \times \fs) \right),
\end{equation} 
where the second equality comes from the definition of $B_i$.

As for \emph{(ii)}, clearly, for all $n \in \N$,
\begin{equation}
  T_o^n \scb_o \times \fs \subseteq \sigma \!\left( \sigma \!\left( 
  \bigcup_{n\in\N} T_o^n \scb_o \right) \times \fs \right),
\end{equation} 
which readily implies the left-to-right inclusion of \emph{(ii)}. On the other hand,
take $B \times \{i\} \in \sigma(\bigcup_{n\in\N} T_o^n \scb_o) \times \fs$. Then
\begin{equation}
  B \times \{i\} \in \sigma \!\left( \bigcup_{n\in\N} (T_o^n \scb_o \times \fs) 
  \right) \subseteq \sigma \!\left( \bigcup_{n\in\N} \sigma(T_o^n \scb_o \times \fs)
  \right).
\end{equation} 
Since the above r.h.s.\ is a $\sigma$-algebra, the right-to-left inclusion of
\emph{(ii)} follows.
\qed

\section{Application to periodic Lorentz gases}
\label{sec-lg}

In this section we use Theorem \ref{thm-main-k}, or rather Corollary 
\ref{cor-main-k}, to show that, under general hypotheses, a periodic Lorentz 
gas is K-mixing \iff it is conservative. This will readily imply that many
2-dimensional periodic Lorentz gases and $d$-dimensional periodic Lorentz 
tubes and slabs are K-mixing. In what follows we use the term 
`Lorentz gas' in its most general meaning of a billiard system defined 
by infinitely many convex scatterers in an unbounded space.

We consider a periodic Lorentz gas in the space $E^d := \R^{d_1} 
\times \mathcal{T}^{d_2}$, where $d_1 \in \Z^+, d_2 \in \N$, $d_1 + d_2 = d$,
and $\mathcal{T}^{d_2}$ is a $d_2$-dimensional torus (not necessarily the 
standard one). So, for $d_1=d$ and $d_2=0$, $E^d$ coincides with $\R^d$. 
See Figs.~\ref{fig1}-\ref{fig3} for examples. 

\newfig{fig1}{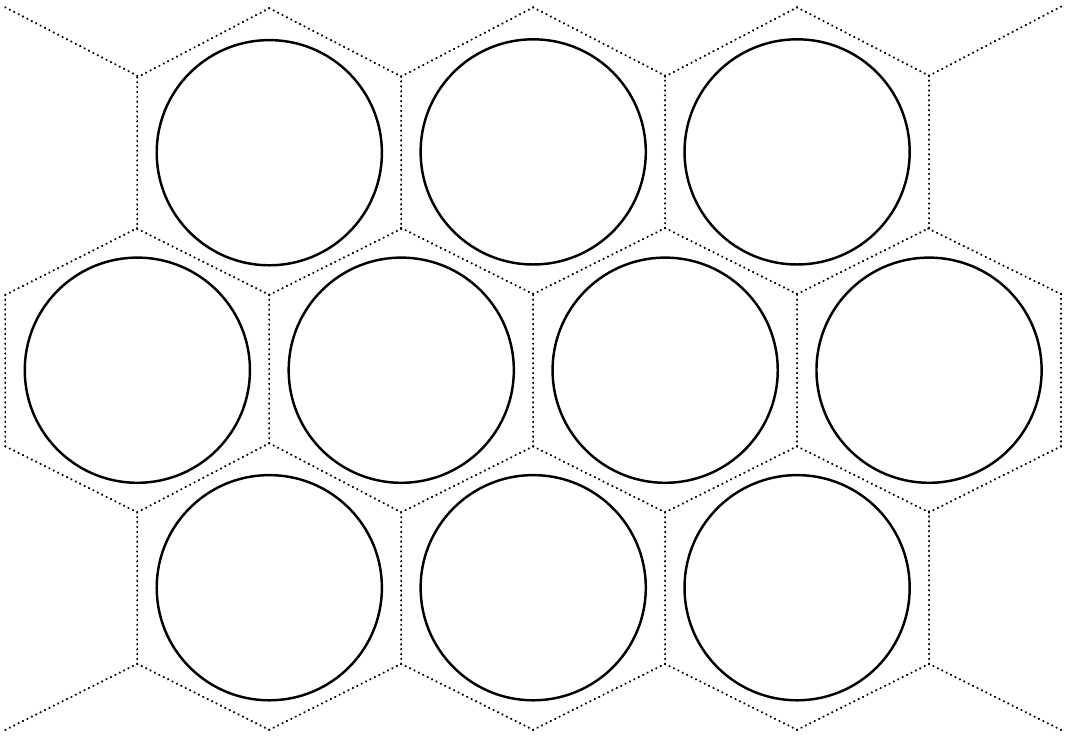}{9cm}{A (proper) periodic Lorentz gas in $\R^2$.}

More in detail, we study the system 
$(\ps, \sca, \mu, T)$ given by a $\Gamma$-extension of the billiard map 
$(\ps_o, \sca_o, \mu_o, T_o)$ of a $d$-dimensional Sinai billiard, where:
\begin{itemize}
\item $\Gamma$ is a co-compact lattice of $\R^{d_1}$, naturally acting on 
  $E^d$. 

\item The Sinai billiard is defined by a finite number of \emph{scatterers} 
  $\mathcal{O}^\ell$, $\ell \in \mathcal{L}$, on $E^d / \Gamma$. The 
  $\mathcal{O}^\ell$ are pairwise disjoint, regular, strictly convex closed sets 
  of $E^d$. (It is not important to prescribe the degree of 
  regularity at this stage, but do observe that $\partial \mathcal{O}^\ell$
  is assumed to be regular and not just piecewise regular. See however 
  Remark \ref{rk-last}.)
  
\item $\ps_o := \bigsqcup_{\ell \in \mathcal{L}} \ps_o^\ell$, where
  \begin{equation}
    \ps_o^\ell := \rset{(q,v) \in \partial \mathcal{O}^\ell \times \R^d} 
    {v \cdot n_q \ge 0, \ |v|=1}
  \end{equation}
  and $n_q$ is the outer unit normal to $\mathcal{O}^\ell$ in $q$. Any such 
  pair $(q,v)$ is also called a \emph{line element} of $ \mathcal{O}^\ell$, 
  and we refer to 
  \begin{equation}
    \Omega_o := (E^d / \Gamma) \setminus \bigsqcup_{\ell \in \mathcal{L}} 
    \mathcal{O}^\ell 
  \end{equation}
  as the \emph{billiard table of the Sinai billiard}.
  
\item $\sca_o$ is the natural Borel $\sigma$-algebra on $\ps_o$ and 
  $\mu_o$ is the usual invariant measure for billiard maps, defined on each 
  $\ps_o^\ell$ by $d \mu_o(q,v) = (v \cdot n_q) dq dv$, where $dq$ is the 
  (hyper)surface element of $\partial \mathcal{O}^\ell$ at $q$, and $dv$ is 
  the (hyper)surface element of $\mathbb{S}^{d-1} := \rset{v \in \R^d} {|v|=1}$ 
  at $v$.
  
\item $T_o: \ps_o \into \ps_o$ is defined so that $T_o(q,v) = (q',v')$ is 
  the line element that the forward trajectory of $(q,v)$ in $\Omega_o$ 
  takes on right after its first collision with some scatterer
  $\mathcal{O}^\ell$.
    
\item $\ps = \ps_o \times \Gamma$ and, for $x = (q,v) \in \ps_o$ and
  $i \in \Gamma$, $\fa(x,i) = i + \psi(x)$ (cf.\ beginning of Section 
  \ref{sec-setup}), where $\psi$ is the \emph{discrete displacement}
  function, defined as follows. We consider the periodic Lorentz gas, 
  that is, the billiard in 
    \begin{equation} \label{omega}
    \Omega := E^d \setminus \bigsqcup_{\ell \in \mathcal{L}}
    \bigsqcup_{i \in \Gamma} \mathcal{O}_i^\ell ,
  \end{equation}
  where $\mathcal{O}_i^\ell$ is the copy of $\mathcal{O}^\ell$ in
  the $i^\mathrm{th}$ copy of the fundamental domain of $\Gamma$ 
  within $E^d$. (In other words, if we tile $E^d$ by means of $\Gamma$
  and draw all the scatterers $\mathcal{O}^\ell$ in the $0^\mathrm{th}$
  tile, then $\mathcal{O}_i^\ell := \mathcal{O}^\ell + i$, a subset of
  the $i^\mathrm{th}$ tile.) So, $\psi(x) = j$ \iff the forward trajectory 
  \emph{in $\Omega$} of $x = (q,v)$, with 
  $q \in \partial \mathcal{O}_0^\ell$, for some $\ell$, has its first collision 
  with a scatterer $\mathcal{O}_j^{\ell'}$, for some $\ell'$. (In other words, 
  the forward trajectory of $(q,v)$ in the periodic Lorentz gas, with
  $q$ interpreted as a point in the $0^\mathrm{th}$ tile, has its first 
  collision in the $j^\mathrm{th}$ tile.) This completely defines $T$. $\sca$ 
  and $\mu$ are then defined out of $\sca_o$ and $\mu_o$ as in Section 
  \ref{sec-setup}. Finally, we refer to $\Omega$ as the \emph{billiard table 
  of the Lorentz gas}, or simply the billiard table.
\end{itemize}

\newfig{fig2}{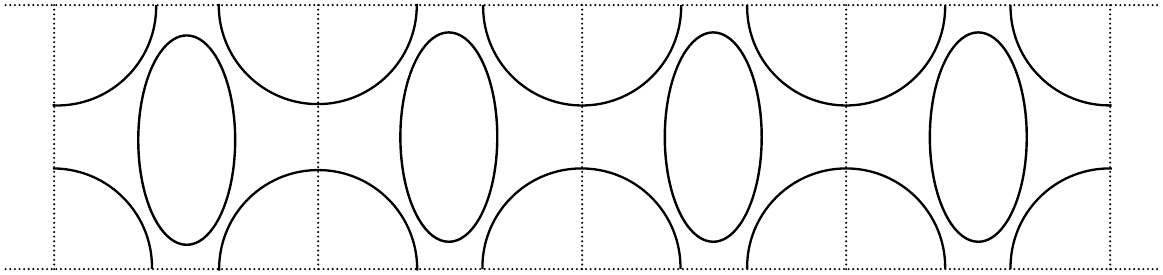}{11cm}{A Lorentz tube in the space $E^2 = \R \times
\mathcal{T}^1$.}

\newfig{fig3}{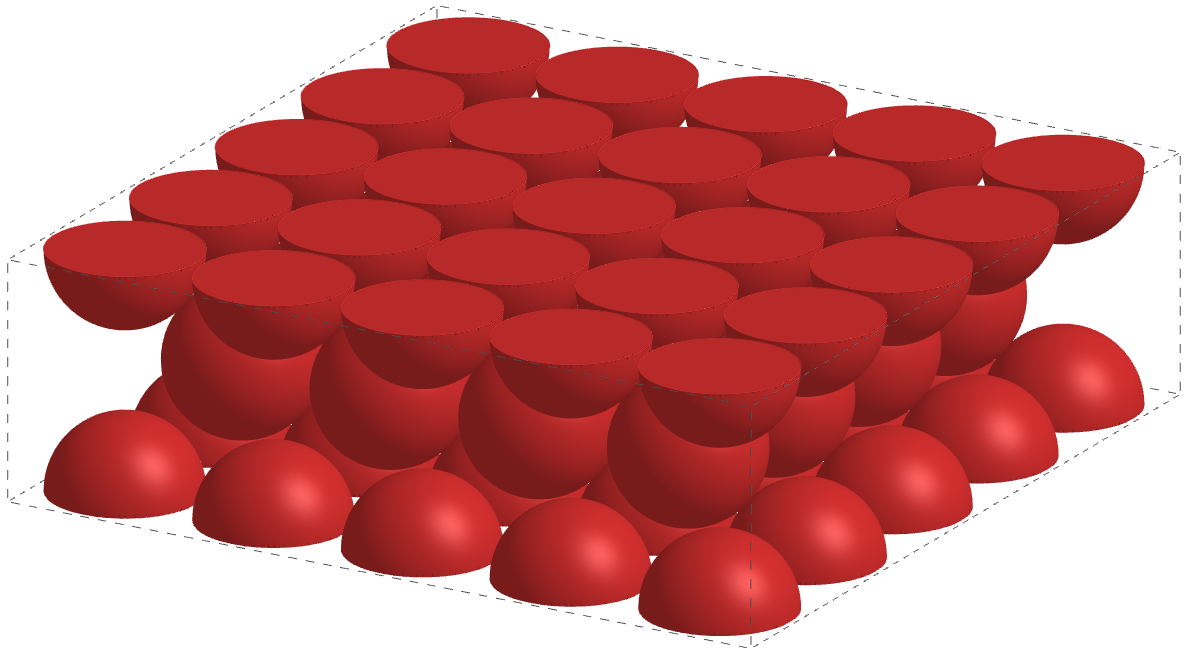}{8cm}{A Lorentz slab in the space $E^3 = \R^2 \times
\mathcal{T}^1$.}

A simple but important observation here is that one could also define
$(\ps, \sca, \mu, T)$ as the scatterer-to-scatterer billiard map for the 
table $\Omega$. In this case $\ps$ would be given as 
$\bigsqcup_{\ell \in \mathcal{L}} \bigsqcup_{i \in \Gamma} \ps_i^\ell$,
where $\ps_i^\ell$ is the set of all line elements for the scatterer 
$\mathcal{O}_i^\ell$; and $T$ would be the first-return map to $\ps$ of
the billiard flow in $\Omega$. Finally, $\mu$ would be the billiard measure 
on the Borel $\sigma$-algebra $\sca$ of $\ps$. This \dsy\ 
is evidently isomorphic to the extension system introduced earlier, so 
much so that we use the same notation for both.

When $d_1=1$ and $d \ge 2$, these Lorentz gases are also referred to as 
\emph{Lorentz tubes}, see Fig.~\ref{fig2}; when $d_1=2$ and $d \ge 3$, they 
are also called \emph{Lorentz slabs}, see Fig.~\ref{fig3}. 

In the interest of generality, the next proposition is formulated in a relatively 
abstract manner, but we will see below that it can be applied easily to many 
examples of periodic Lorentz gases.

\begin{proposition} \label{prop-lg}
  Let $(\ps, \sca, \mu, T)$ be a periodic Lorentz gas corresponding to a
  Sinai billiard $(\ps_o, \sca_o, \mu_o, T_o)$, as described above. Assume
  the following:
  \begin{itemize}
  \item[(a)] $T_o$ has an atomic K-mixing decomposition in the sense of 
    Corollary \ref{cor-main-k}, for some $\sigma$-algebra $\scb_o$
    relative to which the discrete displacement function $\psi$ is measurable.
    
  \item[(b)] $T$ is conservative.
  
  \item[(c)] For all $m \in \Z^+$, $\ell \in \mathcal{L}$ and $i \in \Gamma$, the 
    first-return map of $T^m$ to $\ps_i^\ell$, denoted $(T^m)_{\ps_i^\ell}$,
    is ergodic. (Observe that $(T^m)_{\ps_i^\ell}$ is well-defined 
    $\mu$-almost everywhere because $T^m$ is conservative by (b).)
  \end{itemize}
  Then $T$ is K-mixing relative to $\scb$, the lift of $\scb_o$.
\end{proposition}

\begin{remark}
  In the opposite direction, observe that any \erg\ Lorentz gas is 
  conservative. This is because any measure-preserving, invertible, \erg\ 
  \dsy\ is conservative. In fact, for a measure-preserving system, 
  $\mathcal{C}$ and $\mathcal{D}$ are invariant \cite[Prop.~1.1.6]{a}. If the
  system is also \erg, then either $\ps = \mathcal{C}$ or $\ps = \mathcal{D}$. 
  In the latter case, it is easy to find a wandering set $W$ such that 
  $\bigsqcup_{n\in \Z} T^n W$, which is invariant because $T$ is invertible, 
  has positive but not full measure.
\end{remark}

\proofof{Proposition \ref{prop-lg}}
 By Corollary \ref{cor-main-k}, $\sct_\scb(T)$, the tail $\sigma$-algebra
of $\scb$ w.r.t.\ $T$, is atomic, with generating partition $\scp$. 

We first claim that $\scp$ is coarser than 
$\{ \ps_i^\ell \}_{\ell \in \mathcal{L}, i \in \Gamma}$, that is, every $\ps_i^\ell$ 
is contained in an atom of $\scp$. By \emph{(b)}, every atom of $\scp$ is part of 
an $m$-cycle, for some $m \in \Z^+$.  Thus, if $\ps_i^\ell$ intersects (in a 
positive-measure set) some atom $P$ of order $m$, then $\ps_i^\ell \cap P$ 
is invariant for $(T^m)_{\ps_i^\ell}$. By \emph{(c)}, then, $\ps_i^\ell
\subseteq P$, which proves the claim. Using an imprecise but understandable 
expression, we describe the inclusion $\ps_i^\ell \subseteq P$ by saying that 
the scatterer $\mathcal{O}_i^\ell$ of the Lorentz gas is contained in $P$.

Let us now prove that $T$ is \erg. If it were not, $\ps$ could be split into
two disjoint invariant sets $A$ and $B$ with $\mu(A), \mu(B) > 0$. 
In the imprecise language just introduced, every scatterer of the Lorentz 
gas is contained in either $A$ or $B$. Take two scatterers $\mathcal{O}_1$ 
and $\mathcal{O}_2$ which are contained in $A$ and $B$, respectively, and 
realize the minimum distance between all such pairs of scatterers. The 
existence of $\mathcal{O}_1$ and $\mathcal{O}_2$ is guaranteed by the 
facts that neither $A$ nor $B$ is empty and the configuration of scatterers is 
periodic and locally finite. Now, clearly there exists a segment connecting 
$\mathcal{O}_1$ and $\mathcal{O}_2$ that is tangent to no scatterer of the
Lorentz gas. This is then a portion of a billiard trajectory (called 
\emph{non-singular} in jargon) whose perturbations still connect
$\mathcal{O}_1$ to $\mathcal{O}_2$, implying that $T$ takes a 
positive-measure set from $A$ to $B$, a contradiction. So, $T$ is \erg\ and 
$\ps$ is an $m$-cycle. 

By Theorem \ref{thm-main-k}\emph{(iv)}, it only remains to prove that $m=1$.
Assume by contradiction that $m\ge 2$. Given an atom $P$, set $P_j := T^j P$,
for all $0 \le j \le m-1$. Take a scatterer $\mathcal{O}_0$ contained in $P_0$. 
Varying a line element $(q,v)$ with 
$q \in \partial \mathcal{O}_0$ and $v$ directed outwardly w.r.t.\ $\mathcal{O}_0$, 
we are bound to find a segment $\overline{qq_2}$, with $q_2$ on the boundary
of some scatterer $\mathcal{O}_2$, such that the interior of $\overline{qq_2}$
contains exactly one point of intersection $q_1$ with a scatterer $\mathcal{O}_1$, 
as in Fig.~\ref{fig4}. The segment is then tangent to $\mathcal{O}_1$ 
in $q_1$ and is part of a billiard trajectory which is called \emph{singular}
precisely because of that tangency. The existence of this segment is a 
consequence of the facts that all scatterers are convex (so $\mathcal{O}_0$
cannot be entirely ``surrounded'' by a unique scatterer) and that, for $\mu$-a.e.\ 
$(q,v)$, the half-line defined by $(q,v)$ intersects a scatterer of the Lorentz gas in 
a point different from $q$ (otherwise the base system, the Sinai billiard, would not 
be well-defined). Observe that, since the ambient space $E^d$ may have a toroidal 
factor, $\mathcal{O}_1$ and $\mathcal{O}_2$ may in principle coincide with 
$\mathcal{O}_0$, although we will soon see that this is not the case. In any event, 
the three segments $\overline{qq_1}$, $\overline{qq_2}$, and $\overline{q_1 q_2}$ 
can be independently perturbed into non-singular trajectory segments connecting, 
respectively, $\mathcal{O}_0$ to $\mathcal{O}_1$, $\mathcal{O}_0$ to 
$\mathcal{O}_2$, and $\mathcal{O}_1$ to $\mathcal{O}_2$, cf.\ Fig.~\ref{fig4}. 
Since $\mathcal{O}_0$ is contained in $P_0$, the properties of the $m$-cycle
imply that both $\mathcal{O}_1$ and $\mathcal{O}_2$ are contained in $P_1$ 
and, at the same time, $\mathcal{O}_2$ is contained in $P_2$ (or $P_0$, if 
$m=2$), which is a contradiction.
\qed

\newfig{fig4}{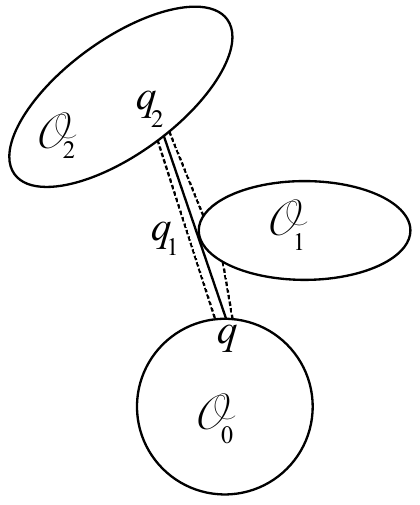}{4.5cm}{Illustration of an argument for the proof of
Proposition \ref{prop-lg}.}

The researcher with a basic familiarity with hyperbolic billiards will see how to
use Proposition \ref{prop-lg} to prove that a given Lorentz gas is K-mixing. 
First of all, the vast majority of Sinai billiards that have been studied (a very 
partial list of references includes \cite{s, bs, ks, sc, bcst, sv, bbt}; see also 
\cite{cm, bl}) are known to be K-mixing w.r.t.\ the $\sigma$-algebra $\scb_o$ 
generated by a 
suitable family of \emph{local stable manifolds (LSMs)}, so ascertaining 
\emph{(a)} likely boils down to a verification in the literature. The lifts to 
$\ps$ of the LSMs for the base system are thus LSMs for the extension and
generate $\scb$ by definition. The first ingredient in the proof of any 
ergodic property for a hyperbolic billiard is the so-called \emph{fundamental 
theorem}, a.k.a.\ \emph{local ergodic theorem}, which in our case reduces to 
the following property: for all $\ell \in \mathcal{L}$, $\mu_o$-a.e.\ pair of points 
$x,y \in \ps_o^\ell$ is connected via a chain of alternating local stable and 
unstable manifolds (LSUMs), with the property that all intersection points 
between successive LSUMs can be chosen in a pre-determined full-measure 
subset of $\ps_o^\ell$. Now, if $T$ is conservative, the LSUMs for $T$ in 
$\ps_i^\ell$ are also LSUMs for $(T^m)_{\ps_i^\ell}$, for all $m \ge1$, 
because, by definition, all points in the same LSM have the same 
\emph{itinerary} in terms of the scatterers hit by the forward trajectories, 
and analogously for LUMs and backward trajectories. Finally, the usual 
\emph{zig-zag argument} for LSUMs (in the terminology of \cite{kss}) shows 
that $(T^m)_{\ps_i^\ell}$ is ergodic, yielding \emph{(c)}.

This method applies easily to all $d$-dimensional Lorentz tubes whose base 
system is K-mixing because of a hyperbolic structure as described above.
In fact, the conservativity of $T$ comes from a known result by Atkinson
\cite{at}, whereby a $\Z$-valued (or $\R$-valued) cocycle over an \erg\ \dsy\ is 
recurrent \iff its displacement function (assumed to be integrable) has zero 
average; see also \cite{cls, sldc}.  As for Lorentz slabs and 2-dimensional 
proper Lorentz gases, independent results by Schmidt \cite{sch} and Conze 
\cite{c} show that recurrence (thus conservativity) is implied by the CLT for 
the displacement function; see also \cite{l2}. So, if the base Sinai billiard is 
K-mixing by way of the standard hyperbolic structure, and satisfies a CLT for 
observables that are locally constant outside its singularity set, the 
corresponding Lorentz slab (or Lorentz gas, in the 2D case) is K-mixing.

\begin{remark} \label{rk-2nd-last}
  For the sake of exposition, we have chosen to discuss only Lorentz 
  gases in $E^d$, which is practically a Euclidean space, but the same ideas 
  and arguments apply to more general ambient spaces, as in the example of 
  Fig.~\ref{fig5}. In fact, $\mathcal{T}^{d_2}$ can be replaced by a polyhedron 
  with $d_1$ pairs of identified parallel facets, with the other facets regarded 
  as hard walls for the billiard dynamics. Proposition \ref{prop-lg} still holds in 
  many cases: one need only check that the scatterers of the Lorentz gas cannot 
  be divided into invariant families (of scatterers) or $m$-cycles with $m \ge 2$ 
  (whose atoms are made up of whole scatterers). The method of application of 
  Proposition \ref{prop-lg} is then the same as described above, showing for 
  example that the Lorentz tube of Fig.~\ref{fig5} is K-mixing.
\end{remark}

\newfig{fig5}{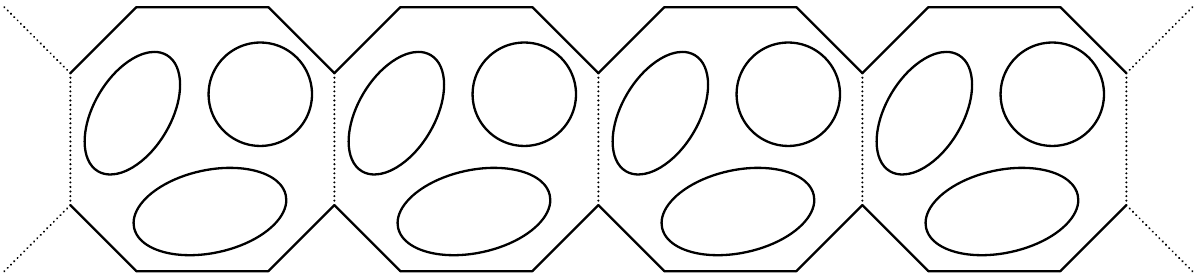}{11cm}{A Lorentz tube with hard walls (cf.\ 
Remark \ref{rk-2nd-last}).}
  
\begin{remark} \label{rk-last}
  One could also extend the above technique to certain Lorentz gases whose
  scatterers are only piecewise regular. An analogue of Proposition \ref{prop-lg} 
  can be proved in some cases, with the difference that, in place of 
  $\ps_i^\ell$, one uses $\ps_i^{\ell,j}$, the part of phase space associated to 
  the $j^\mathrm{th}$ regular component of $\partial \mathcal{O}_i^\ell$. 
  (The proof might be cumbersome in this case, and depends on which regular 
  boundary components of other scatterers are ``seen'' by each regular boundary 
  component of each scatterer.) It is not hard to see that the whole procedure
  works for the Lorentz tube of Fig.~\ref{fig6}, which is therefore K-mixing.
\end{remark}

\newfig{fig6}{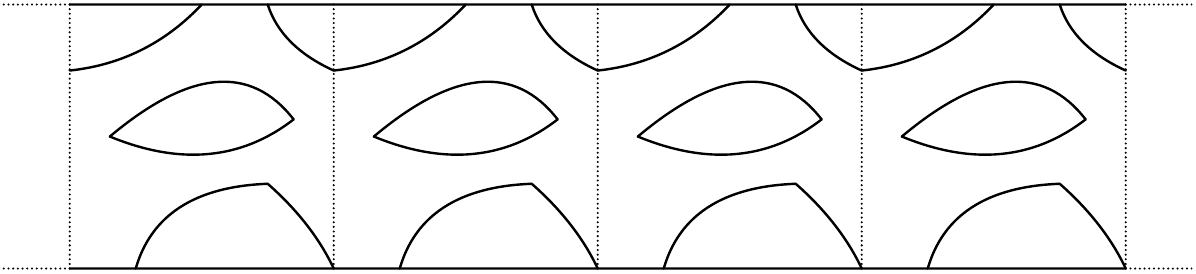}{11cm}{A Lorentz tube with (hard walls and) piecewise 
regular scatterers (cf.\ Remark \ref{rk-last}).}

\footnotesize


\begin{thebibliography}{BCST}

\bibitem[A]{a} \book{J.~Aaronson} {An introduction to infinite ergodic theory} 
  {Mathematical Surveys and Monographs, 50. American Mathematical 
  Society, Providence, RI, 1997}

\bibitem[At]{at} \article{G.~Atkinson} {Recurrence for co-cycles and
  random walks} {J. London.  Math. Soc. (2) \vol{13} (1976),
  486--488}

\bibitem[BBT]{bbt} \article{P.~Bachurin, P.~B\'alint and I.~P.~T\'oth}
  {Local ergodicity for systems with growth properties including 
  multi-dimensional dispersing billiards} {Israel J. Math. \vol{167} (2008), 
  155--175}

\bibitem[BCST]{bcst} \article{P.~B\'alint, N.~Chernov, D.~Sz\'asz and
  I.~P.~T\'oth} {Multi-dimensional semi-dispersing billiards:
  singularities and the fundamental theorem} {Ann. Henri Poincar\'e
  \vol{3} (2002), no.~3, 451--482}

\bibitem[BK]{bk} \article{T.~Bogensch\"utz and Z.~S.~Kowalski} {Exactness 
  of skew products with expanding fibre maps} {Studia Math. \vol{120} (1996), 
  no.~2, 159--168}
   
\bibitem[BS]{bs} \article{L.~A.~Bunimovich and Ya.~G.~Sinai} {On a 
  fundamental theorem in the theory of dispersing billiards} {Math. USSR-Sb.
   \vol{19} (1973), 407--423. Translated from Mat. Sbornik N.S. \vol{90(132)} 
   (1973), 415--431}

\bibitem[BL]{bl} \article{L.~Bussolari and M.~Lenci} {Hyperbolic billiards with 
  nearly flat focusing boundaries. I} {Physica D \vol{237} (2008), no.~18, 
  2272--2281}

\bibitem[CM]{cm} \article{N.~Chernov and R.~Markarian} {Chaotic
  billiards} {Mathematical Surveys and Monographs, 127. American
  Mathematical Society, Providence, RI, 2006}

\bibitem[C]{c} \article{J.-P.~Conze} {Sur un crit\`ere de r\'ecurrence en 
  dimension 2 pour les marches stationnaires, applications} {Ergodic 
  Theory Dynam. Systems \vol{19} (1999), no.~5, 1233--1245}

\bibitem[CLS]{cls} \article{G.~Cristadoro, M.~Lenci and M.~Seri}
  {Recurrence for quenched random Lorentz tubes} {Chaos \vol{20}
  (2010), 023115. Errata corrige in Chaos \vol{20} (2010), 049903}

\bibitem[KS]{ks} \book{A.~Katok and J.-M.~Strelcyn (in collaboration 
  with F.~Ledrappier and F.~Przytycki)} {Invariant manifolds, entropy 
  and billiards; smooth maps with singularities} {Lectures Notes in
  Mahematics 1222, Springer-Verlag, Berlin-New York, 1986}
  
\bibitem[K]{k} \article{Z.~S.~Kowalski} {The exactness of generalized 
  skew products} {Osaka J. Math. \vol{30} (1993), no.~1, 57--61}

\bibitem[KSS]{kss} \article{A.~Kr\'amli, N.~Sim\'anyi and D.~Sz\'asz}
  {A ``transversal'' fundamental theorem for semidispersing billiards}
  {Comm. Math. Phys. \vol{129} (1990), no.~3, 535--560. Errata corrige
  in Comm. Math. Phys. \vol{138} (1991), no.~1, 207--208}

\bibitem[L1]{l1} \article{M.~Lenci} {Aperiodic Lorentz gas: recurrence
  and ergodicity} {Ergodic Theory Dynam. Systems \vol{23} (2003),
  no.~3, 869--883}

\bibitem[L2]{l2} \article{M.~Lenci} {Typicality of recurrence for
  Lorentz gases} {Ergodic Theory Dynam. Systems \vol{26} (2006),
  no.~3, 799--820}
  
\bibitem[L3]{l3} \article{M.~Lenci} {Exactness, K-property and infinite 
  mixing} {Publ. Mat. Urug. \vol{14} (2013), 159--170} 
  
\bibitem[L4]{l4} \article{M.~Lenci} {Uniformly expanding Markov maps of 
  the real line: exactness and infinite mixing} {Discrete Contin. Dyn. Syst. 
  \vol{37} (2017), no.~7, 3867--3903} 
  
\bibitem[P]{p} \article{Ya.~B.~Pesin} {Characteristic Lyapunov exponents 
  and smooth ergodic theory} {Russian Math. Surveys \vol{32} (1977), no.~4, 
  55--114. Translated from Uspehi Mat. Nauk \vol{32} (1977), no.~4, 55--112}
  
\bibitem[R]{r} \article{V.~A.~Rohlin} {On the fundamental ideas of 
  measure theory} {Amer. Math. Soc. Translation \vol{1952} (1952), 
  no.~71, 55 pp. Translated from Mat. Sbornik N.S. \vol{25(67)} (1949),
  107--150}

\bibitem[Sch]{sch} \article{K.~Schmidt} {On joint recurrence} {C. R. Acad.
  Sci. Paris S\'er. I Math. \vol{327} (1998), no.~9, 837--842}

\bibitem[SLDC]{sldc} \article{M.~Seri, M.~Lenci, M.~Degli Esposti and 
  G.~Cristadoro} {Recurrence and higher ergodic properties for quenched 
  random Lorentz tubes in dimension bigger than two} {J. Stat. Phys. 
  \vol{144} (2011), no.~1, 124-138}
  
\bibitem[S]{s} \article{Ya.~G.~Sinai} {Dynamical systems with elastic 
  reflections} {Russ. Math. Surveys \vol{25}, 137--189 (1970). Translated 
  from Uspehi Mat. Nauk \vol{25} (1970) no.~2, 141--192}
  
\bibitem[SC]{sc} \article{Ya.~G.~Sinai and N.~I.~Chernov} {Ergodic 
  properties of certain systems of two-dimensional discs and 
  three-dimensional balls} {Russ. Math. Surveys \vol{42} (1987), no.~3, 
  181--207 (1987). Translated from Uspekhi Mat. Nauk \vol{42} (1987), 
  no.~3, 153--174}
  
\bibitem[SV]{sv} \article{D.~Sz\'asz and T.~Varj\'u} {Limit laws and 
  recurrence for the planar Lorentz process with infinite horizon.} {J. Stat. 
  Phys. \vol{129} (2007), no.~1, 59--80}
  
\end{thebibliography}
\end{document}